\magnification=\magstep1
\input amstex
\documentstyle{amsppt}
\voffset=-3pc
\loadbold
\loadmsbm
\loadeufm
\UseAMSsymbols
\baselineskip=12pt
\parskip=6pt
\def\var{\varepsilon}
\def\bC{\Bbb C}
\def\bN{\Bbb N}
\def\bR{\Bbb R}

\def\bZ{\Bbb Z}
\def\cO{\Cal O}

\def\cU{\Cal U}
\def\cV{\Cal V}

\def\oDelta{\overline\Delta}

\NoBlackBoxes
\topmatter
\title Analytic Continuation in Mapping Spaces\endtitle
\address Department of Mathematics,
Purdue University,
West Lafayette, IN\ 47907-1395\endaddress
\author
L\'aszl\'o Lempert\endauthor\footnote""{Research partially supported by NSF
grant DMS0203072\hfill\break}
\abstract
We consider a Stein manifold $M$ of dimension $\geq 2$ and a compact subset $K\subset M$ such that $M'=M\backslash K$ is connected.
Let $S$ be a compact differential manifold, and let $M_S$, resp.~$M'_S$ stand for the complex manifold of maps $S\to M$, resp.~$S\to M'$, of some specified regularity that are homotopic to constant.
We prove that any holomorphic function on $M'_S$ continues analytically to $M_S$ (perhaps as a multivalued function).
\endabstract
\endtopmatter
\document

\subhead 1.\ Introduction\endsubhead

A fruitful idea in geometry is to study a space $M$ by considering maps from a fixed space $S$ to $M$, and relating properties of the 
resulting mapping space to those of $M$.
For example, homotopy groups of $M$ are defined in terms of the space of continuous maps $S^i\to M$.
If $M$ is a complex manifold, so will be its mapping space(s); motivated by a theorem of Hartogs and Serre on analytic continuation in Stein manifolds, in this paper we shall investigate analytic continuation in mapping spaces.
Let $M$ be a Stein manifold of dimension $m\geq 2$, and $K\subset M$ a compact subset such that $M'=M\backslash K$ is connected.
The Hartogs--Serre theorem then says that any holomorphic function on $M'$ continues to a holomorphic function on $M$, see [H, Se].
(Hartogs was the first to discover this phenomenon and prove its 
first instances; he also formulated a general theorem on analytic 
continuation in domains in $\bC^n$.
However, he seems to have worked with multivalued functions, and even with this interpretation his proof is not quite convincing.
The general result for Stein manifolds was first given by Serre.)

We shall consider a mapping space version of this result, and prove an analogous theorem, but there will be a difference.
Let $S$ be a compact differential manifold, possibly with boundary, and $A\subset S$ closed.
Fix furthermore a regularity class $C^k\ (k=0,1,2,\ldots)$ or Sobolev $W^{k,p}\ (k=1,2,\ldots,1\leq p<\infty,kp>\dim S)$\footnote"$^{*}$"{The Sobolev embedding theorem thus will guarantee that elements of $W^{k,p}(S)$ are represented by continuous functions.}, and also a point $0\in M'$.
Denote by $M_{S,A}$ the space of maps $x\colon (S,A)\to (M,0)$, of the given regularity, that are homotopic to the constant map; and define $M'_{S,A}$ similarly.
It turns out that $M_{S,A}$ has a natural structure of a complex Banach manifold---in fact this is true for any complex manifold $M$, see [L2],---and $M'_{S,A}\subset M_{S,A}$ is open.
Our main result is

\proclaim{Theorem 1.1}In addition to the assumptions and notation above, 
if $A\neq\emptyset$
assume also that there is a proper holomorphic map $\Delta\to M$ of the unit disc $\Delta\subset\bC$ that passes through $0\in M$ and avoids $K$.
If $Y$ is a complex Banach space, then any holomorphic function $f\colon M'_{S,A}\to Y$ can be continued analytically along any curve in $M_{S,A}$, starting at $x_0\equiv 0$.
\endproclaim

The problem of analytic continuation in mapping spaces was
first considered in 1972 by Greenfield, who obtained a variant
of the above theorem in the special case when 
$(M,K)=(\Bbb C^n,\{0\})$ and $A=\emptyset$, see [Gr]. 
He focused on the space of 
continuous maps, and allowed $S$ to be a compact metric space.

In Section 3 we shall recall notions surrounding analytic continuation. In Theorem 1.1
the condition on the proper map $\Delta\to M$ is not very restrictive.
Loosely speaking, it is satisfied if $0$ is sufficiently far from 
$K$.---When $S\backslash A$ is a singleton, the pair 
$(M_{S,A}, M'_{S,A})$
is biholomorphic to the pair $(M,M')$, so that the Hartogs--Serre 
theorem applies:
the analytic continuation along curves gives rise to a single valued holomorphic function $M_{S,A}\to Y$, extending $f$.
This is also true in the case considered in [Gr]. However,
in general the situation is more complicated.
Whether $f$ has a single valued holomorphic extension 
$M_{S,A}\to Y$ depends on whether its analytic continuation along closed curves is single valued.
One is led to consider closed curves through $x_0$, along which any holomorphic $f\colon M'_{S,A}\to\bC$ has a single valued analytic continuation.
The homotopy classes of such curves constitute a subgroup $\Gamma=\Gamma(M_{S,A},M'_{S,A})\subset\pi_1(M_{S,A},x_0)$, that we call the monodromy group.
If $\Gamma=\pi_1(M_{S,A},x_0)$ then any holomorphic $f\colon M'_{S,A}\to\bC$ extends holomorphically to $M_{S,A}$ (and the same will be true for Banach space valued $f$, as a straightforward application of the Banach--Hahn theorem shows).
We have been able to identify $\Gamma$ for mapping spaces based on most, but not all, regularity classes.

\proclaim{Theorem 1.2}Suppose our mapping spaces are defined using a Hilbertian Sobolev regularity class $W^{k,2}$ ($2k>\dim S$), or any other regularity class $C^l$ or $W^{l,p}$ contained in it.
If $A\neq\emptyset$, also suppose that there is a proper holomorphic map $\Delta\to M$ through $0\in M$ that avoids $K$.
Then the monodromy
group $\Gamma$ is the image of the homomorphism $\pi_1 (M'_{S,A},x_0)\to\pi_1 (M_{S,A},x_0)$ induced by the inclusion $M'_{S,A}\to M_{S,A}$.
\endproclaim

The same is true for certain other regularity classes, but we cannot deal with mapping spaces defined by continuous maps.
The reason is that the proof of Theorem 1.2 depends on certain recent results in infinite dimensional complex geometry that are currently available only in Banach spaces which have an unconditional basis; yet 
(for positive dimensional $S$) the space $C(S)$ does not have an
unconditional basis.

For a while I was in the mistaken belief that $\Gamma=\pi_1(M_{S,A},x_0)$, and, as a result, $\pi_1(M'_{S,A},x_0)\to\pi_1(M_{S,A},x_0)$ is onto.
However, Gompf pointed out that this is not always so (and I am grateful to Gong, who told me about Gompf's example).
Theorem 1.2 is still of interest.
It provides a new instance of how topological and complex analytic properties determine each other in the theory of Stein spaces.

\subhead 2.\ Background\endsubhead

In this section we collect a few general results about complex Banach manifolds, in particular, about mapping spaces.
For the notion of holomorphic maps in (always complex) Banach spaces, complex Banach manifolds, and holomorphic maps between them, we refer to [L1, Section 2].
In this paper we shall only consider what we termed in [L1] rectifiable complex Banach manifolds, i.e.~those that are locally biholomorphic to open subsets of Banach spaces (and we shall drop the adjective ``rectifiable'').
We denote the set of holomorphic maps between complex Banach manifolds by $\cO(M;N)$.

Consider an open $\Omega\subset\bC^n$ and a complex submanifold $N\subset\Omega$.
For the sake of simplicity, we assume that $N$ contains the origin $0\in\bC^n$.
Suppose that $\rho\colon\Omega\to N$ is a holomorphic retraction, i.e.~$\rho|N=\text{ id}_N$.
Let $S$ be a compact smooth manifold, possibly with boundary; and $A\subset S$ closed.
As in the Introduction, fix a regularity class $C^0$, $C^k$, or $W^{k,p}$, $k=1,2,\ldots$, $kp>\dim S$, and define the mapping space $N_{S,A}$ to consist of those maps $(S,A)\to (N,0)$ of the given regularity that are homotopic to the constant map $x_0\equiv 0$.
In particular, $\bC^n_{S,A}$ is a Banach space and $\Omega_{S,A}\subset\bC^n_{S,A}$ is open.

\proclaim{Proposition 2.1}$N_{S,A}\subset\Omega_{S,A}$ is a direct complex submanifold, and the map
$$
\Omega_{S,A}\ni x\mapsto\rho\circ x\in N_{S,A}\tag2.1
$$
is a holomorphic retraction on it.
\endproclaim

Recall that a closed subset $P$ of a complex Banach manifold $Q$ is a complex submanifold if for every $x\in P$ there are a neighborhood $\cU\subset Q$ of $x$, a Banach space $X$, a closed subspace $Y\subset X$, an open $\cV\subset X$, and a biholomorphic map $\Phi\colon\cU\to\cV$ such that $\Phi(P\cap\cU)=Y\cap\cV$.
If $Y$ has a closed complement in $X$, then $P$ is called a direct submanifold.

\demo{Proof}Let $g\colon\bC^n\times N\to T \bC^n|N$ denote the standard trivialization: $g(\zeta,z)\in T_z\bC^n$ is the velocity vector of the curve $r(t)=z+t\zeta$ at $t=0$.
The kernel of $d\rho$ along $N$ defines a normal bundle to $N$, i.e.~a holomorphic subbundle $E\to N$ of $T\bC^n|N$ such that $T\bC^n|N=E\oplus TN$.
With the projections
$$
\pi_E\colon T\bC^n|N\to E\qquad\text{ and }\qquad\pi_T\colon T\bC^n|N\to TN
$$
define a holomorphic map $f\colon N\times\bC^n\to T\bC^n|N$ by
$$
f(\zeta,z)=\pi_T g(\zeta,z)+\pi_E g(\zeta+z-\rho(\zeta+z),z).
$$
Thus $f$ is a fiberwise map, and $f(0,\cdot)$ is the zero section of $T\bC^n|N$.
The difference
$$
f(\zeta,z)-g(\zeta,z)=\pi_E g(\zeta+z-\rho(\zeta+z),z)-\pi_E g(\zeta,z)
$$
vanishes to second order along $\{0\}\times N$, as one checks separately along tangent vectors to $\bC^n\times \{z\}$ parallel to $T_z N$, resp.~$E_z$.
Hence $df=dg$ is an isomorphism at points of $\{0\}\times N$, and the inverse function theorem (plus an exercise in paracompactness, see e.g.~[DG, or L3, p.~84]) implies that $f$ is a fiberwise biholomorphism between neighborhoods $U\subset\bC^n\times N$ of $\{0\}\times N$ and $V\subset T\bC^n|N$ of the zero section.
We choose $U$ so that $\zeta+z\in\Omega$ if $(\zeta,z)\in U$.

We are now ready to show that $N_{S,A}\subset\Omega_{S,A}$ is a direct submanifold.
Indeed, it is clearly closed.
If $x\in N_{S,A}$, define $X$ to be the Banach space of those $C^k$, resp.~$W^{k,p}$ sections of the induced bundle $x^*T\bC^n\to S$ that vanish on $A$.
Let $Y,Z\subset X$ be the spaces of sections valued in $x^*TN$, resp.~$x^* E$.
Thus $X=Y\oplus Z$.
Let furthermore
$$
\gathered
\cU=\{y\in\Omega_{S,A}\ \colon \ (y(s)-x(s),\ x(s))\in U\text{ for all }s\in S\},\\
\cV=\{\sigma\in X\colon\sigma(s)\in V\text{ for all }s\in S\}.
\endgathered
$$
One easily checks that
$$
\Phi\colon \cU\ni y\mapsto f(y-x,x)\in\cV
$$
is a biholomorphism and $\Phi (N_{S,A}\cap\cU)=Y\cap\cV$; which means that $N_{S,A}$ is a direct submanifold.

The second part of the proposition is obvious.
\enddemo

Further down we shall need an analogous result about a space of 
mappings of another type, this time into a manifold that 
itself may be infinite dimensional.
Consider a Banach space $(X,\|\ \|)$, an open $\Omega\subset X$, and a holomorphic retraction $\rho\colon\Omega\to N$ on a complex submanifold $N\subset\Omega$.
Let $\overline\Delta$ be a finite dimensional compact complex manifold, possibly with boundary, and interior $\Delta$ (the only case we shall later need is $\oDelta\subset\bC$ the closed unit disc).
Continuous maps $x\colon\oDelta\to X$ that are holomorphic on $\Delta$ constitute a Banach space $X_\Delta$, with norm given by $\sup_{s\in\oDelta} \|x(s)\|$, $x\in X_\Delta$.
Consider the space $\Omega_\Delta$, resp.~$N_\Delta$ of those $x\in X_\Delta$ that map into $\Omega$, resp.~$N$.
Then $\Omega_\Delta\subset X_\Delta$ is open.

\proclaim{Proposition 2.2}$N_\Delta\subset\Omega_\Delta$ is a direct complex submanifold and the map
$$
\Omega_\Delta\ni x\mapsto\rho\circ x\in N_\Delta
$$
is a holomorphic retraction on it.
\endproclaim

This is proved the same way as Proposition 2.1.

Finally, we shall need an approximation result for curves in manifolds.
Let again $X$ be a Banach space, $\Omega\subset X$ open, $\rho\colon\Omega\to N$ a holomorphic retraction on a submanifold.

\proclaim{Proposition 2.3}If $I\subset\bR$ is a compact interval, any continuous map $f\colon I\to N$ can be uniformly approximated by analytic maps $g\colon I\to N$.
\endproclaim

That $g$ is analytic means that $I\subset\bC$ has a neighborhood $U$ to which $g$ extends as a holomorphic map with values in $N$.

\demo{Proof}By a slight extension of Weierstrass' approximation theorem, there are polynomials $P_\nu\colon\bC\to X$ such that $P_\nu\to f$ uniformly on $I$.
(Most proofs of the original theorem give this vector valued version as well, e.g.~the one in [Ru, pp.~159--160] does.)
The maps $g_\nu=\rho\circ P_\nu|I$ are then analytic for $\nu > \nu_0$ and converge to $f$ uniformly.
\enddemo

\subhead 3.\ Generalities on analytic continuation\endsubhead

One ingredient in the proof of Theorem 1.1 is the method of sliding discs, an idea that goes back to Hartogs.
We shall employ it in the following setting.
Let $X,Y$ be Banach spaces, $O\subset X$ open, $r\colon O\to N$ a holomorphic retraction on a complex submanifold $N\subset O$, and $N'\subset N$ open. In what follows, by a curve we 
shall mean a continuous map of $[0,1]$ in some space.

\definition{Definition 3.1}We say that an $f\in\cO(N';Y)$ can be continued analytically along a curve $[0,1]\ni t\mapsto x_t\in N$, starting at $x_0\in N'$, if for each $t\in [0,1]$ a neighborhood $V_t\subset N$ of $x_t$ and $f_t\in\cO (V_t; Y)$ can be given so that $f_0=f$ near $x_0$, and with some $\delta>0$, $f_t=f_s$ on $V_t\cap V_s$ whenever $|t-s|<\delta$.
\enddefinition

Equivalently, if $\cO^Y\to N$ denotes the sheaf of $Y$--valued holomorphic germs, one can define analytic continuation as a continuous map $\varphi\colon [0,1]\to\cO^Y$ that covers $t\mapsto x_t$, and $\varphi(0)$ is the germ of $f$ at $x_0$.
Analytic continuation is unique in the sense that if $\{V'_t,f'_t\}$ define another analytic continuation of $f$ along the same curve,
then $f_t=f'_t$ near $x_t$.
The identity theorem implies that if the curve $t\mapsto x_t$ is deformed, with the endpoints kept fixed, and $f$ can be continued along each of the deformed curves, then the germ of the continuation $f_1$ at $x_1$ will not change.

Let $\Delta\subset\bC$ be the unit disc, and $N_\Delta$ the complex Banach manifold of continuous maps $\xi\colon\oDelta\to N$ that are holomorphic on $\Delta$, as in Proposition 2.2.
Elements of $N_\Delta$ will be called (analytic) discs.
The distance between two discs $\xi,\eta$ is
$$
\|\xi-\eta\|_\Delta=\sup_{\sigma\in\overline\Delta} \|\xi(\sigma)-\eta(\sigma)\|,
$$
where $\|\ \|$ stands for the norm on $X$.

\definition{Definition 3.2}A curve $[0,1]\ni t\mapsto\xi_t\in N_\Delta$ 
is called regular, if $\xi_t(\partial\Delta)\subset N'$ for all 
$t\in [0,1]$, and $\xi_0(\overline\Delta)\subset N'$.
It is a regular lift of a curve $[0,1]\ni t\mapsto x_t\in N$ if, in 
addition, $\xi_t (0)=x_t$ for $t\in [0,1]$.
\enddefinition

\proclaim{Lemma 3.3}For any regular curve $t\mapsto\xi_t\in N_\Delta$ there is an $\varepsilon>0$ with the following property.
If $[0,1]\ni t\mapsto\eta_t\in N_\Delta$ and $[0,1]\ni t\mapsto\zeta_t\in N_\Delta$ are curves that satisfy $\|\xi_t-\eta_t\|_\Delta<\var$, $\|\xi_t-\zeta_t\|_\Delta<\var\quad (0\leq t\leq 1)$, and $\eta_u(0)=\zeta_u(0)$ with some $u\in [0,1]$, then $t\mapsto\eta_t,\zeta_t$ are regular, and
for every $f\in\cO(N';Y)$
$$
\int_0^1 f(\eta_u (e^{2\pi i\tau}))d\tau=\int_0^1 f(\zeta_u (e^{2\pi i\tau}))d\tau.\tag3.1
$$
\endproclaim

\demo{Proof}Let
$$
K=\xi_0 (\oDelta)\cup\bigcup_{0\leq t\leq 1}\xi_t(\partial\Delta)\subset N',
$$
a compact set.
Choose $\varepsilon > 0$ so that the distance between $\bigcup_{0\leq t\leq 1}\xi_t(\oDelta)$ and $\partial O$ is $> 3\varepsilon$; and for any $x\in K$ and $y\in X$ with $\|y\| < 3\varepsilon$ we have $r(x+y)\subset N'$.
This guarantees that $t\mapsto\eta_t,\zeta_t$ are regular.
To show (3.1), approximate $\eta_t,\zeta_t$ by analytic curves $t\mapsto\eta_t^*$, $t\mapsto\zeta^*_t$ (cf. Proposition 2.3),
and define the analytic curve
$$
[0,1]\ni t\mapsto\omega_t=r(\eta_t^*-\eta_t^* (0)+\zeta_t^* (0))\in N_\Delta.
$$
If the approximation is close enough, $t\mapsto\zeta^*_t,\omega_t$ will
be regular curves. As $\zeta_t^*(\oDelta)$,
$\omega_t(\oDelta)\subset N'$ for $t$ in some neighborhood of 0,
for such $t$ $f$ is holomorphic on $\zeta_t^*(\oDelta)$, $\omega_t(\oDelta)$, whence
$$
\int_0^1 f(\zeta_t^* (e^{2\pi i\tau}))d\tau=f(\zeta_t^*(0))=f(\omega_t (0))=\int_0^1 f(\omega_t (e^{2\pi i\tau}))d\tau.\tag3.2
$$
Since both integrals here depend analytically on $t\in [0,1]$, it follows that they agree for $t=u$ as well.
Letting $\eta_t^*$ and $\zeta_t^*$ go to $\eta_t$, resp $\zeta_t$, we obtain (3.1). 
\enddemo

\proclaim{Lemma 3.4 (of sliding discs)}Any $f\in\cO (N';Y)$ can be analytically continued along a curve $[0,1]\ni t\mapsto x_t\in N$ starting at $x_0\in N'$, that has a regular lift.
\endproclaim

\demo{Proof}Let $\xi_t$ be a regular lift of $x_t$.
Fix $\var>0$ as in Lemma 3.3.
There is a $\delta_1>0$ such that whenever $y\in N$ is at distance $<\delta_1$ to some $x_t$, there is a regular curve $[0,1]\ni u\mapsto\eta_u\in N_\Delta$, at distance $<\var/2$ to $\xi_u$, such that $\eta_t(0)=y$ (for example, we can take
$$
\eta_u=r(\xi_u - x_t + y),\qquad 0\leq u\leq 1).\tag3.3
$$
Next choose $\delta>0$ so that
$$
\|\xi_t-\xi_s\|_\Delta < \var/2,\qquad \text{if}\quad |t-s| < \delta.
\tag3.4
$$

To continue $f$ analytically, let
$$
V_t=\{ y\in N\colon \| x_t - y\| < \delta_1\},
$$
and for $y\in V_t$ define $f_t(y)$ in the following way.
Construct a regular curve $[0,1]\ni u\mapsto\eta_u\in N_\Delta$ at distance $<\var/2$ to the curve $u\mapsto\xi_u$, such that $\eta_t(0)=y$, and put
$$
f_t (y)=\int_0^1 f (\eta_t (e^{2\pi i\tau}))d\tau\in Y.\tag3.5
$$
By Lemma 3.3 $f_t(y)$ does not depend on which curve $\eta_u$ we choose.
Now $f_t$ is holomorphic.
Indeed, if we use the lift (3.3), the integral in (3.5) is manifestly a holomorphic function of $y\in V_t$.

Next suppose that we use the lift (3.3) to compute $f_0(y)$.
Since $f$ is holomorphic on $\eta_0(\oDelta)$, (3.5) gives $f_0(y)=f(y)$.
Finally, suppose $|t-s| < \delta$, and let $y\in V_t\cap V_s$.
We compute $f_t(y)$ by choosing a curve $u\mapsto \eta_u\in N_\Delta$ as described above.
To compute $f_s(y)$ we use the curve $u\mapsto\zeta_u=\eta_{u+t-s}\in N_\Delta$.
(To make sense of this, we extend $\eta_u$ constant for $u\leq 0$ and for $u\geq 1$.)
Then with $\xi_u$ similarly extended
$$
\|\zeta_u-\xi_u\|_\Delta\leq \|\zeta_u-\xi_{u+t-s}\|_\Delta + \|\xi_{u+t-s}-\xi_u\|_\Delta<\var,
$$
in view of (3.4).
Also $\zeta_s(0)=y$.
The upshot is that $f_s(y)$ can be computed as
$$
f_s(y)=\int_0^1 f(\zeta_s (e^{2\pi i\tau}))d\tau=\int_0^1 f(\eta_t (e^{2\pi i\tau}))d\tau=f_t (y),
$$
which completes the proof.
\enddemo

The Lemma of sliding discs allows one analytically to continue holomorphic functions on $N'$ along any curve in $N$ that starts at a fixed 
$x_0\in N'$, provided one can show all such curves have regular lifts.
However, analytic continuation along curves does not guarantee that an 
$f\in\cO (N';Y)$ can be continued to
a holomorphic function on $N$, because continuations along curves will, in general, give rise to a multivalued function.
The rest of the section will be devoted to this issue.

We shall say that an analytic continuation $\{(V_t,f_t)\colon 0\leq t\leq 1\}$ along a closed curve $t\mapsto x_t$ is single valued if $f_0=f_1$ in a neighborhood of $x_0=x_1$.
Assuming that $f$ can be analytically continued along any such curve, whether the analytic continuation is single valued depends only on the homotopy class of the curve.
We fix $x_0\in N'$, and define the monodromy group 
$\Gamma=\Gamma(N,N')\subset\pi_1 (N,x_0)$ to consist of those 
homotopy classes along which analytic continuation of any 
$f\in\cO(N'; \bC)$ is single valued.
In fact, for any Banach space $Y$, analytic continuation of any
$f\in\cO(N'; Y)$ along a homotopy class in $\Gamma$ will be single valued
(since for any linear form $\varphi:Y\to\bC$ the analytic continuation
of $\varphi f$ is such).

\proclaim{Proposition 3.5}Suppose that every $f\in\cO(N';\bC)$ can be analytically continued along every curve $[0,1]\ni t\mapsto x_t\in N$ starting at a fixed $x_0\in N'$.
Then the monodromy group $\Gamma$ contains the image of 
$\pi_1(N', x_0)\to\pi_1 (N,x_0)$.
\endproclaim

\demo{Proof}Any homotopy class in the image is represented by a closed curve $t\mapsto x_t\in N'$, and analytic continuation along the latter is obviously single valued.
\enddemo

Often much more can be said:

\proclaim{Theorem 3.6}Suppose that $X$ has an unconditional basis and 
$O\subset X$ is pseudoconvex (see e.g.~[LP 1.4, 1.5]).
If for every Banach space $Y$ every $f\in\cO (N',Y)$ can be analytically
continued along every curve in $N$, starting at $x_0\in N'$, then the 
monodromy group $\Gamma$ is the image of $\pi_1(N',x_0)\to\pi_1 (N,x_0)$.
\endproclaim

To verify this we need the following result, generalizing K.~Stein's theorem that coverings of Stein manifolds are Stein [St].

\proclaim{Theorem 3.7}Let $X,O$ be as in Theorem 3.6, $P\subset O$ a direct complex submanifold, and $\pi\colon Q\to P$ a holomorphic covering.
Then $Q$ can be embedded in some Banach space $Y$ as a direct complex submanifold.
If $\pi$ has countably many sheets, then $Y$ can be taken to have an unconditional basis.
\endproclaim

\demo{Proof}By Zerhusen's theorem in $[Z]$, $O$ can be embedded as a direct complex submanifold in a Banach space $X_1$ with unconditional basis.
Since this embedding sends $P$ to a direct complex submanifold of $X_1$, we can assume right away that $O=X$ and so $P$ is a submanifold of $X$.
We also assume $P$ is connected.

First we construct a holomorphic Banach bundle $E\to P$ out of $Q$.
Its fibers $E_x$ are the Banach spaces $l^1(\pi^{-1} x),\ x\in P$.
If $\pi|\pi^{-1} U$ is trivial for some open $U\subset P$, i.e.~it is isomorphic to the projection $\Sigma\times U\to U$ with some set $\Sigma$, then this isomorphism induces a bijection
$$
\coprod_{x\in U} E_x\to l^1 (\Sigma)\times U.
$$
We endow the set theoretical vector bundle $E=\coprod_{x\in P} E_x\to P$ with the holomorphic Banach bundle structure for which these bijections are local trivializations.

There is another holomorphic Banach bundle $E'\to P$ such that 
$E\oplus E'=F$ is trivial, see [Pa1, Theorem 1.3b and Pa2, Theorem 11.1a].
The former reference is about holomorphic Banach bundles over 
pseudoconvex open sets $\Omega\subset X$; in the latter the bundles are over submanifolds.
This more general result follows from the former and from the existence of holomorphic neighborhood retractions, [LP, Theorem 12.3 or Pa2, Theorem 1.4f].
In fact, the fibers of ($E'$ and) $F$ are countable sums of the fibers $E_x$, completed with respect to the $l^1$ norm.
If $\pi$ has countably many sheets, it follows that the fibers of $F$ also have unconditional bases.

There is an obvious holomorphic embedding 
$$
d\colon Q\to E\subset F
$$
that associates with $q\in Q$ the element $d(q)\in l^1(\pi^{-1}\pi q)$ that is 1 on $q$ and 0 everywhere else.
It is easy to check that $d(Q)$ is a direct complex submanifold of $F$.
Finally, $F$ can be extended to a trivial bundle $F'\to X$, whose total space is biholomorphic to a Banach space $Y$; then $d(Q)$ becomes a direct complex submanifold of $F'$, and so $Q$ has the required embedding into $Y$.
Since $Y\approx F_x\times X$, it also follows that $Y$ has an unconditional basis, provided $\pi$ has countably many sheets.
\enddemo

\demo{Proof of Theorem 3.6}Let the image of $\pi_1(N',x_0)\to\pi_1 (N,x_0)$ be $G$; it determines a covering $\pi\colon Q\to N$.
As a topological space $Q$ is the quotient of the space of paths in $N$ starting at $x_0$, two paths $\gamma_1,\gamma_2$ being equivalent if they have the same endpoint, and the loop obtained by concatenating $\gamma_1$ and $-\gamma_2$ has homotopy class in $G$.
The projection $\pi$ associates with any $q\in Q$ the endpoint of (a representative of) $q$.
Note that a closed curve $t\mapsto x_t$ in $N$ lifts to a closed curve in $Q$ precisely when its homotopy class is in $G$.
Since $\pi$ is a local homeomorphism, the complex structure of $N$ lifts to $Q$.
Theorem 3.7 shows that $Q$ can be considered a submanifold of a Banach space $Y$.
Observe that $\pi$ has a holomorphic section $f$ over $N'$, associating with $x\in N'$ the class $q=f(x)\in\pi^{-1} x$ of a curve from $x_0$ to $x$ that runs in $N'$.

Suppose a closed curve $[0,1]\ni t\mapsto x_t\in N$ represents a homotopy class $\gamma\in\Gamma$.
This implies that along it $f\colon N'\to Q\subset Y$ has single valued analytic continuation $\{V_t,f_t\}$, $f_t\in\cO(V_t,Y)$.
In fact, $f_t$ will map some neighborhood of $x_t\in V_t$ into $Q$, and $\pi\circ f_t=$ id will hold there for all $t\in [0,1]$.
This can be verified by showing that those $t$ for which this is true constitute a closed and open subset of $[0,1]$.
But then $t\mapsto f_t(x_t)\in Q$ is a closed curve that lifts $t\mapsto x_t$; and as we have noted, this means that $\gamma\in G$.
Therefore $\Gamma\subset G$, and, in light of Proposition 3.5, in fact $\Gamma=G$.
\enddemo

Our final result also belongs to a section about analytic continuation, even though it is not applicable in the context of Theorem 1.1.

\proclaim{Theorem 3.8}Suppose that $H^1(N,\cO^Y)=0$, and that there is a 
$g\in\cO(N;\bC)$ vanishing at a fixed $x_0\in N'$, such that 
$\inf_{N\backslash N'} |g|>0$.
If every $f\in\cO (N';Y)$ can be analytically continued along every 
curve starting at $x_0$, then the analytic continuations form a single 
valued continuation $\tilde f\in\cO(N;Y)$ of $f$.
\endproclaim

\demo{Proof}All we have to show is that the analytic continuation of $f$ along any closed curve $[0,1]\ni t\mapsto x_t\in N$ is single valued.
Suppose the continuation is given by a family $V_t$ of open sets and $f_t\in\cO(V_t;Y)$.
We can assume that $[0,1]$ is partitioned in finitely many intervals $I_j$ such that for each $j$, $V_t$ does not depend on $t\in I_j$; and that $V_t\supset x(\overline I_j)$ for $t\in I_j$.
In this case the analytic continuation of $f$ along any curve $t\mapsto y_t$, sufficiently close to $t\mapsto x_t$, can also be given by the family $(V_t,f_t)$.
We choose this nearby curve $t\mapsto y_t$ to start at $x_0$, go inside $V_0$ to some $y$, and eventually end at $y$.
For points $y$ in some neighborhood of $x_0$ such a curve will exist.
If, in addition, $|g(y)| < \inf_{N\backslash N'} |g|$, we shall show that $f_0 (y)=f_1 (y)$.

To this end let $N''=\{x\in N\colon g(x)\neq g(y)\}$, and 
$$
h=f/(g-g(y))\in\cO (N'\cap N'';Y).
$$
Since $N'\cup N''=N$, by our cohomological assumption we can write $h=h'+h''$ with $h'\in\cO(N';Y)$ and $h''\in\cO(N'';Y)$.
In other words
$$
f=(g-g(y)) h'+ (g-g(y)) h''\qquad\text{on } N'\cap N''.
$$
This formula defines a holomorphic extension $f''$ of the last term to $N'\cup N''=N$.
On the other hand, continue $h'$ analytically along the curve 
$t\mapsto y_t$.
The identity $f=(g-g(y)) h'+f''$ will be preserved in the process, and we obtain
$$
f_1(y)=(g(y)-g(y))h'_1 (y)+f''(y)=f''(y)=f_0(y).
$$
Therefore the analytic continuation of $f$ along the curve $t\mapsto x_t$ is indeed single valued, which proves the theorem.  
\enddemo

\subhead 4.\ Analytic discs in Stein Manifolds\endsubhead

To apply the results of Section 3 in mapping spaces, we will study the space of analytic discs in Stein manifolds.
Consider a Stein manifold $M$ of dimension $\geq 2$,
properly embedded in $\bC^n$.
According to Docquier and Grauert [DG], $M$ has a pseudoconvex neighborhood 
$\Omega\subset\bC^n$ that admits a holomorphic retraction 
$\rho\colon\Omega\to M$.
Let $u\colon M\to\bR$ be a smooth, strongly plurisubharmonic exhaustion function, and for real numbers $a<b$ define
$$
\aligned
M(a)&=\{z\in M\colon u(z) < a\},\\
M(a,b)&=\{z\in M\colon a< u(z) < b\},\\
M[a,b]&=\{z\in M\colon a\leq u(z)\leq b\}.
\endaligned
$$
We fix $a<b$ so that $u$ has no critical value on $[a,b]$.
As in Section 3, $\Delta$ will denote the unit disc in $\bC$, and $M_\Delta$ the complex manifold of continuous maps $\oDelta\to M$ that are holomorphic on $\Delta$.
Finally, let
$$
\gathered
D(a,b)=\{x\in M_\Delta\ \colon \ x(0)\in M(a),\ x(\partial\Delta)\subset M(a,b)\},\qquad\text{and}\\
\pi\colon D(a,b)\ni x\mapsto x(0)\in M(a).
\endgathered
$$
By the maximum principle $x(\oDelta)\subset M(b)$ for $x\in D(a,b)$.~---That we can apply the method of sliding discs in mapping spaces $M_{S,A}$ depends on the following result.

\proclaim{Theorem 4.1}$\pi\colon D(a,b)\to M(a)$ is a Serre fibration.
\endproclaim

This means that for any cube $Q=[0,1]^d$, maps $Q\to M(a)$ have the homotopy lifting property:\ given continuous maps
$$
h\colon [0,1]\times Q\to M(a)\quad\text{ and }\quad 
H_0\colon Q\to D(a,b)
$$
such that $\pi\circ H_0=h(0,\cdot)$, there is a continuous map
$H\colon [0,1]\times Q\to D(a,b)$ such that $H(0,\cdot)=H_0$ and $\pi\circ H=h$.

It is not even obvious that $\pi$ of Theorem 4.1 is surjective, but this, and certain more precise results about discs in Stein manifolds were 
proved by Forstneri\v c and Globevnik in [FG,G].
Our proof was largely inspired by their work.
This proof will be given in Section 6;
here and in the next section we present a few auxiliary results.

We denote by $E\to M[a,b]$ the smooth complex subbundle of $TM[a,b]$, of corank 1, consisting of vectors in the kernel of $\partial u$.
Let $\nu\colon M[a,b]\to TM[a,b]$ be the zero section.

\proclaim{Lemma 4.2}There are a fiberwise convex open neighborhood $V\subset E$ of $\nu(M[a,b])$ and a smooth map $\varphi\colon V\to M$, holomorphic on the fibers $V_z=V\cap E_z$, such that $\varphi\circ\nu=\text{ id}_{M[a,b]}$, and for any $z\in M[a,b]$ and real line $l\subset E_z$ through $\nu(z)$ the function $u\circ\varphi|l\cap V_z$ has a unique critical point at $\nu(z)$, a minimum.
\endproclaim

\demo{Proof}For brevity, we write $N=M[a,b]$ and $v=u\circ\rho\colon\Omega\to\bR$.
Using coordinates in $\bC^n$, for $\zeta,z\in N$ let 
$$
q(\zeta,z)=2\sum^n_{i=1}\ {\partial v(z)\over\partial z_i}\ (\zeta_i-z_i)+\sum^n_{i,j=1}\ {\partial^2 v(z)\over\partial z_i\partial z_j}\ (\zeta_i-z_i)(\zeta_j-z_j).
$$
Since $\partial_\zeta q (\zeta,z)=\partial v(\zeta)\neq 0$ when $\zeta=z$, the implicit function theorem gives a neighborhood $U\subset N\times N$ of the diagonal such that $L=\{(\zeta,z)\in U\colon q(\zeta,z)=0\}$ is a smooth submanifold of real codimension 2.
Also, $L_z=\{\zeta\in N\colon (\zeta,z)\in L\}$ is a complex submanifold of $U_z=\{\zeta\in N\colon (\zeta,z)\in U\}$, of codimension 1.
Clearly, $T_z L_z=E_z$.
By Taylor expansion, as $N\ni \zeta\to z$ 
$$
u(\zeta)=v(\zeta)=v(z)+\text{Re }q(\zeta,z)+\sum^n_{i,j=1}\ {\partial^2 v(z)\over\partial z_i\partial\overline z_j}\ (\zeta_i-z_i)(\overline\zeta_j-\overline z_j)+O(|\zeta-z|^3)
$$
Since $(\partial^2v/\partial z_i\partial\overline z_j)$ is 
positive definite on $T_z N$, it follows that $u|L_z$ has a 
nondegenerate local minimum at $z$.

Now choose a hermitian metric on $T\bC^n$,
and let $\pi_E\colon T\bC^n|N\to E$ be the orthogonal projection.
As in the proof of Proposition 2.1 we consider the trivialization 
$g\colon\bC^n\times N\to T\bC^n| N$, but this time define a smooth map
$$
f\colon L\ni (\zeta,z)\mapsto\pi_E g(\zeta-z,z)\in E.
$$
One checks that $f$ maps $L_z$ holomorphically to $E_z$, $f(z,z)=\nu(z)$, and---since $f-g$ vanishes to second order
at $(z,z)$---$df$ defines an isomorphism between $T_{(z,z)} L_z$ and 
$T_{\nu(z)} E_z$ for $z\in N$.
It follows that $f$ is a diffeomorphism between a neighborhood of $\{(z,z)\colon z\in N\}$ in $L$ and a neighborhood $V\subset E$ of $\nu(N)$.
Composing its inverse $V\to L\subset N\times N$ with the projection $N\times N\to N$ on the first factor we obtain a map $\varphi$ with $\varphi\circ\nu=\text{ id}_{M[a,b]}$, and the restrictions $u\circ\varphi|V_z$ have nondegenerate local minima at $\nu(z)$.
If $V$ is now suitably shrunk, then $u\circ\varphi|l\cap V_z$ will have a unique critical point at $\nu(z)$ for every real line $l\subset E_z$, as needed.
\enddemo

\proclaim{Lemma 4.3}There is a $\delta_0>0$ with the following property.
Let $t_0<t_1$ be real numbers, and $Q$ a compact cube.
Suppose $\kappa\colon Q\times\oDelta\to M$ is a continuous map,  that
is holomorphic on each slice $\{q\}\times\Delta$, $q\in Q$, and
maps $Q\times\partial\Delta$ in $M(a,b)$; 
$\alpha\colon Q\times\partial\Delta\to (0,\delta_0]$ is continuous; and 
$\eta$ is a positive number.
Then there is a continuous map 
$\lambda\colon [t_0,t_1]\times Q\times\oDelta\to M$, holomorphic on each slice $\{(t,q)\}\times\Delta$, such that
\roster
\item"(i)" 
$\kappa(q,s)=\lambda(t_0,q,s)\text{ for }(q,s)\in Q\times\oDelta$;
\item"(ii)"
$\kappa(q,0)=\lambda(t,q,0)\text{ for }(t,q)\in [t_0,t_1]\times Q$;
\item"(iii)" 
$\alpha(q,s)-\eta < u(\lambda (t_1,q,s))-u(\kappa (q,s)) 
< \alpha (q,s)+\eta\text{ for }(q,s)\in Q\times\partial\Delta$;
\item"(iv)" 
$-\eta < u(\lambda (t,q,s))-u(\kappa(q,s))
<\alpha (q,s)+\eta\text{ for }(t,q,s)\in 
[t_0,t_1]\times Q\times\partial\Delta$.
\endroster
\endproclaim

This is a variant of [G, Lemma 5.2] and of [FG, Main Lemma].
The punchline, formulae (4.5), (4.6), is the same as there.
The underlying idea made its first appearance in several complex variables 
in Poletsky's papers, see e.g.~[Po, p.~168]; but, according to Poletsky, 
in one variable it had been used earlier by Lavrentiev.

\demo{Proof}Let $E,V,\varphi$ be as in Lemma 4.2, and assume $t_0=0,\ t_1=1$.
Choose $\delta_0>0$ so that for every $z\in M[a,b]$ 
$$
\{w\in V_z\colon u(\varphi(w))< u(z)+\delta_0\}\subset\subset V_z.
$$

Set $\kappa_1=\kappa|Q\times\partial\Delta$.
Since $Q\times\partial\Delta$ is homotopically equivalent to
the circle, and over the circle every complex vector bundle of
finite rank is trivial,
the induced bundle $\kappa_1^* E$ is trivial, and has a 
nonvanishing section.
This section gives rise to a continuous map 
$g\colon Q\times\partial\Delta\to E$ that covers $\kappa_1$, 
and avoids the zero section.
It follows that for $(q,s)\in Q\times\partial\Delta$
$$
\Gamma_{qs}=\{\sigma\in\bC\ \colon \ u(\varphi (\sigma g(q,s)))=
u(\kappa (q,s))+\alpha (q,s)\}\tag4.1
$$
is a smooth Jordan curve, that depends continuously on $(q,s)$.
Let $f_{qs}$ denote the biholomorphic map of $\Delta$ on the inside
$$
\{\sigma\in\bC \ \colon \ u(\kappa(q,s))\le 
u(\varphi(\sigma g(q,s)))< u(\kappa (q,s))+\alpha (q,s)\}\tag4.2
$$
of $\Gamma_{qs}$, normalized by $f_{qs}(0)=0$, $f'_{qs}(0)>0$.
These maps extend to homeomorphisms of $\oDelta$, and the extended maps, also denoted $f_{qs}$, depend continuously on $q,s$ (see [C, Ra]).
Consider the continuous map
$$
\psi\colon Q\times\partial\Delta\times\oDelta\ni (q,s,\sigma)
\mapsto\varphi\bigl(f_{qs}(\sigma g(q,s))\bigr)\in M\subset\bC^n,\tag4.3
$$
holomorphic in $\sigma\in\Delta$.
It can be uniformly approximated by continuous $\bC^n$ valued maps of form
$$
\chi(q,s,\sigma)=\kappa (q,s)+\sum_{|i|< J}\sum^J_{j=1}\ a_{ij} (q) s^i\sigma^j.
$$
Indeed, viewing $\psi$ as a function of arg $s$, arg $\sigma$, and expanding it in a Fourier series, the (iterated) Ces\`aro means of the partial sums will be of the above form, and converge uniformly to $\psi$, according to Fej\'er's theorem.
Choose $\chi$ so that it maps into $\Omega$ and
$$
|u\circ\psi - u\circ\rho\circ\chi| < \eta\qquad\text{ on }Q\times\partial\Delta\times\oDelta.\tag4.4
$$

Observe that
$$
\chi(q,s,(ts)^J)=
\kappa (q,s)+\sum_{|i| < J}\sum^J_{j=1} a_{ij} (q) s^{i+jJ} t^{jJ};\tag4.5
$$
the right hand side is holomorphic in $s\in\Delta$, and maps 
$s\in\partial\Delta$ in $\Omega$.
Since when $t=0$, it maps any $s\in\oDelta$ in $\Omega$, the 
Kontinuit\"atsatz (see e.g.~[K, Theorem 3.3.5]) implies it maps
any $(t,q,s)\in[0,1]\times Q\times\oDelta$ in $\Omega$.
In particular,
$$
\lambda(t,q,s)=
\rho\bigl(\kappa (q,s)+
\sum_{|i|< J}\sum^J_{j=1} a_{ij}(q) s^{i+j J} t^{jJ}\bigr)
\tag4.6
$$
is holomorphic in $s\in\Delta$ and satisfies (i), (ii).
It also follows from (4.1),\dots, (4.6) that (iii), (iv) hold, q.e.d.
\enddemo

\subhead 5.\ Harmonic measure\endsubhead

In this section we discuss a few facts about harmonic measure in the complex plane, that will be used in the proof of Theorem 4.1.

\definition{Definition 5.1}Let $G\subset\bC$ be a bounded domain, $0\in G$, and $\Gamma\subset\partial G$ a Borel set.
Given a positive number $\delta$, we say that the harmonic measure of $\Gamma$ relative to $G$ is $>\delta$ if there is a harmonic function $\omega\colon G\to\bR$ such that
$$
\omega(0)=\delta,\ \qquad\text{and}\qquad\overline{\lim\limits_{s\to s_0}} \omega(s) <
\cases 1,&\text{if $s_0\in\Gamma$}\\ 0,&\text{if $s_0\in\partial G\backslash\Gamma$}.\endcases\tag5.1
$$
\enddefinition

Loosely speaking, harmonic measure is a holomorphically invariant way of measuring length on boundaries surrounding 0.
More formally, if $\psi\colon\overline G'\to\overline G$ is continuous, $\psi|G'\colon G'\to G$ is biholomorphic, and $\psi(0)=0$, then relative to $G$ a set $\Gamma\subset\partial G$ has harmonic measure $>\delta$ precisely when $\psi^{-1}(\Gamma)$ has, relative to $G'$.~---We shall write
$$
\Delta(\var)=\{s\in\bC\colon |s|<\var\},\tag5.2
$$
with $\Delta=\Delta(1)$ as before.

\proclaim{Proposition 5.2}If a relatively open $\Gamma\subset\partial\Delta$ has harmonic measure $>\delta$ (relative to $\Delta$), then there is a complex polynomial $\theta$ such that
$$
\theta(0)=\delta,\quad\text{Re }\theta < \cases 1&
\text{on $\overline\Delta$}\\
0&\text{on $\partial\Delta\backslash\Gamma$.}\endcases
$$
\endproclaim

\demo{Proof}Let $\omega$ be as in Definition 5.1, and extend it to 
$\partial\Delta$ by the $\overline\lim$ in (5.1).
This extension, also denoted $\omega$, is upper semicontinuous; on 
the other hand, the characteristic function 
$\chi_\Gamma\colon\partial\Delta\to\bR$ is lower semicontinuous.
Hence by Hahn's insertion theorem there is a continuous function, even 
a (trigonometric) polynomial $\omega_1\colon\partial\Delta\to\bR$, 
such that $\omega<\omega_1<\chi_\Gamma$ on $\partial\Delta$.
Let $\theta_1$ be a complex polynomial on $\bC$ with $\text{Re }\theta_1|\partial\Delta=\omega_1$.
By the mean value theorem $\delta=\omega(0)<\text{ Re }\theta_1 (0)$; therefore $\theta=\theta_1-\theta_1(0)+\delta$ is the required polynomial.
\enddemo

\proclaim{Lemma 5.3}For every $\var\in (0,1)$ there is a $\delta>0$ with the following property.
Let $\sigma_0\in\partial\Delta$, $\eta<1$, and
$$
\Sigma_j=\{s\in\Delta\backslash\{0\}\colon (j-1)\var <\text{ arg }s/\sigma_0 < j\var\},\qquad j=0,1.
$$
Let furthermore $G\subset\Delta$ be a simply connected domain,
$$
G\supset\Delta (\var)\cup\Sigma_0\cup \{s\in\Sigma_1\colon |s|>\eta\}.
$$
If $\partial G$ intersects the open segment $0\sigma_0$, then relative to $G$ the harmonic measure of $\Sigma_1\cap\partial G$ is $>\delta$.
\endproclaim

\demo{Proof}Consider the open set $U=\Delta(\var)\cup \{s\in\bC\ \colon \ 0 <\text{ arg }s<\var\}$ and its image $V$ under the biholomorphic map $g(s)=1-\sqrt{1-s}$; we use the branch for which $g(0)=0$.
Both $U$ and $V$ are unbounded and connected, hence so is $(V\backslash\{0\})\cup (\bC\backslash\overline{\Delta(3/\var)})$.
We can apply Runge's approximation theorem on the complement $\{0\}\cup (\overline{\Delta(3/\var)}\backslash V)$, to obtain a polynomial $\theta$ on $\bC$ such that $\theta(0)=1$ and Re $\theta<0$ on $\overline{\Delta(3/\var)}\backslash V$.
We claim that any $\delta < 1/\sup_{\Delta(3/\var)}\text{ Re }\theta$ will do.

To verify this, let $\sigma_0,G$ be as in the Lemma.
By symmetry we can assume $\sigma_0=1$, so that there is a $\sigma\in (\var,1)\cap\partial G$.
Define
$$
g_\sigma(s)=g(s/\sigma)=1-\sqrt{1-s/\sigma},\qquad s\in U;
$$
then $g_\sigma(U)\supset V$.
On the other hand, consider the single valued branch $f(s)$ of $1-\sqrt{1-s/\sigma}$ on $G$, chosen so that $f(0)=0$.
Then $G'=f(G)\subset\Delta(3/\var)$.
The inverse of $f$ is given by $f^{-1}(t)=\sigma t(2-t)$, hence it extends to a continuous function $\psi\colon\overline G'\to\overline G$.
We will show that the harmonic measure of $\Gamma'=\psi^{-1}(\Sigma_1\cap\partial G)$ relative to $G'$ is $>\delta$; this will then imply the Lemma since harmonic measure is holomorphically invariant.

The point is that $V\cap\partial G'\subset\Gamma'$, or, equivalently, that $\psi(V\cap\partial G')\subset\Sigma_1\cap\partial G$.
Indeed, if $t\in V\cap\partial G'$, then $\psi(t)\in\partial G$.
Further, $\psi(t)=g_\sigma^{-1}(t)\in U$, since the functions
$\psi$ and $g_\sigma^{-1}$ agree where both are defined.
Now $U\cap\partial G$ consists of $\Sigma_1\cap\partial G$ and of an arc of $\partial\Delta$.
Suppose $\psi(t)=g_\sigma^{-1}(t)=s$ were on this arc.
Consider a closed curve $\gamma\subset\Delta$, starting at $s$, going to $e^{-i\var/2}$ along the shorter arc $I$ of $\partial\Delta$, then straight to 0, and finally straight to $s$.
Now $f$ can be extended across $I$, and along $\gamma$, $f$ so extended and $g$ are analytic continuations of one another.
They would map $\gamma$ to a closed curve if $f(s)=g_\sigma(s)(=t)$.
However, this is impossible because $\gamma$ surrounds the branch point $\sigma$ precisely once.
We conclude that $\psi(t)\in\Sigma_1\cap\partial G$, i.e., $\psi(V\cap\partial G')\subset\Sigma_1\cap\partial G_1$, and $V\cap\partial G'\subset\Gamma'$ as claimed.

It now follows that the harmonic function $\omega=\delta$ Re $\theta$ satisfies $\omega(0)=\delta$ and
$$
\omega<\cases 1&\text{on $\overline{G'}\subset\overline{\Delta(3/\var)}$}\\
0&\text{on $\partial G'\backslash\Gamma'\subset\overline{\Delta(3/\var)}\backslash V$.}\endcases
$$
Therefore the harmonic measure of $\Gamma'$ relative to $G'$ is $>\delta$, and so is the harmonic measure of $\Sigma_1\cap\partial G$ relative to $G$.
\enddemo

\subhead 6.\ The Proof of Theorem 4.1\endsubhead

We start with a simple topological result:

\proclaim{Lemma 6.1}Let $T=P_0P_1 P_2$ be a closed triangle in the plane, covered by relatively open subsets $G_1,G_2$.
If $G_i\supset P_0 P_i$ for $i=1,2$, then $P_0$ can be connected with the side $P_1 P_2$ by a simple polygon running in $G_1\cap G_2$.
\endproclaim

\demo{Proof}There are closed $F_i\subset G_i$ containing the sides $P_0 P_i$ in their relative interior that still cover $T$; in fact, we can arrange that they are bounded by finitely many simple closed polygons.
Consider the boundary component of $F_1$ that contains $P_0 P_1$; removing $P_0 P_1$ from it, we obtain a simple polygon $\Pi$ running from $P_0$ to $P_1$, but otherwise avoiding the side $P_0 P_1$.
Let $P$ be the first point on $\Pi$ that is also on $P_1P_2$, and $\Gamma$ the part of $\Pi$ between $P_0$ and $P$.
Obviously, $\Gamma\subset F_1\subset G_1$; but also $\Gamma\subset F_2\subset G_2$.
Indeed, let $R\in\Gamma$.
If $R\in P_0 P_2$ then of course $R\in F_2$.
Otherwise arbitrarily close to $R$ there are points in $T\backslash F_1\subset F_2$, which again implies $R\in F_2$.
Hence $\Gamma$ is the polygon sought.
\enddemo

Now we resume the notation and assumptions of Theorem 4.1.
To prove the theorem, we have to consider a cube $Q=[0,1]^d\subset\bR^d$, a homotopy $h\colon [0,1]\times Q\to M(a)$, and a lift $H_0\colon Q\to D(a,b)$ of $h(0,\cdot)$; and we have to construct a lift $H\colon [0,1]\times Q\to D(a,b)$ of $h$ so that $H(0,\cdot)=H_0$.
We fix $h$, and in Propositions 6.2, 6.3, 6.4 all constants that occur may depend on $h$.

\proclaim{Proposition 6.2}There is a $\delta_1>0$ such that if 
$x\in D(a,b)$ and $x(0)\in h([0,1]\times Q)$, then 
$x(\Delta(\delta_1))\subset M(a)$, cf.~(5.2).
\endproclaim

\demo{Proof}If $M(a,b)\subset\bC^n$ is contained in a ball of radius $r$ and $x\in D(a,b)$, then $|x(s)-x(0)|\leq 2r|s|$ by Schwarz's lemma.
Therefore $\delta_1$ will do if the distance between $h([0,1]\times Q)$ and $M\backslash M(a)$ is $> 2r\delta_1$.
\enddemo

\proclaim{Proposition 6.3}There is a $\delta_2>0$ with the following property.
Suppose that, with some $\tau\in [0,1)$ and $a< a' < c < b' < b$, 
$h(\tau,\cdot)$ has been lifted to $K_\tau\colon Q\to D(a',b')$.
Then there are $\overline\tau\in (\tau,1)$, arbitrarily close to 
$\tau$, and a lift $K\colon [\tau,\overline\tau]\times Q\to D(a',b')$ of 
$h|[\tau,\overline\tau]\times Q$, such that $K(\tau,\cdot)=K_\tau$, and 
for any $q\in Q$, $x=K(\overline\tau,q)$ maps a set 
$\Gamma\subset\partial\Delta$ of harmonic measure $>\delta_2$ (relative 
to $\Delta$) into $M(c)$.
\endproclaim

\demo{Proof}Fix $\delta_1$ as in Proposition 6.2, and with $\var=\min(\delta_1, 2^{-d})$, let $\delta$ be as in Lemma 5.3.
We shall show that $\delta_2=\delta$ is a possible choice.
To construct the lift $K$ we shall compose $K_\tau(q)$ with certain 
holomorphic maps $f_{tq}\colon\Delta\to\Delta$, $q\in Q$, 
$t\in [\tau,\overline\tau]$, to obtain a map like $K$, except it will 
not lift $h|[\tau,\overline\tau]\times Q$; but this map can be perturbed 
to a lift, provided $\overline\tau$ is sufficiently close to $\tau$.

The maps $f_{tq}$ will in fact map $\Delta$ biholomorphically on
certain domains $G_{tq}\subset\Delta$, and will be obtained as follows.
Let $q\in Q$ and $x=K_\tau (q)$.
Consider the sectors
$$
\Sigma_j=\{s\in\Delta\backslash \{0\}\colon j 2^{-d} < \text{ arg }s < (j+1) 2^{-d}\},\quad j=0,1,\ldots,3\cdot 2^d-1,
$$
with vertices $0, s_j=\exp (ij 2^{-d})$ and $s_{j+1}$.
Topologically $\overline{\Sigma}_{3j+1}\cup\overline\Sigma_{3j+2}$ is a triangle $T$ with vertices $0, s_{3j+1}, s_{3j+3}$.
We cover it by relatively open sets
$$
\align
G_1&=\{s\in T\colon\quad x(s)\in M(a',b')\},\\
G_2&=\{s\in T\colon\quad x(s)\in M(c)\text{ or }s\not\in\overline\Sigma_{3j+1}\}.
\endalign
$$
We are in the situation of Lemma 6.1:\ the side $s_{3j+1} s_{3j+3}$ is covered by $G_1$, the side $0s_{3j+3}$ by $G_2$.
We conclude that there is a simple polygon $\Gamma_{jq}\subset G_1\cap G_2$ connecting $s_{3j+3}$ with a point on the side $0s_{3j+1}$.
We can arrange that $\Gamma_{jq}\backslash \{s_{3j+3}\}$ is in $\Delta$.
Thus
$$
x(s)\in\cases M(a',b'),&\text{if $s\in\Gamma_{jq}$}\\ M(c),&\text{if $s\in\Gamma_{jq}\cap\overline\Sigma_{3j+1}$}.\endcases\tag6.1
$$
Given $\overline\tau\in (\tau,1)$ and an $\eta>0$, parametrize each $\Gamma_{jq}$ with $t\in [0,\overline\tau-\tau]$, starting at $s_{3j+2}$; let $\Gamma_{jq}(t)\subset\Gamma_{jq}$ denote the arc corresponding to $[0,t]$, and $F_{qt}(t)$ the closed neighborhood of $\Gamma_{jq}(t)$ in $\overline\Sigma_{3j+1}\cup\overline\Sigma_{3j+2}$, of radius $\min(t,\eta)$.
We write $F_{jq}=F_{jq}(\overline\tau-\tau)$, so that $F_{jq}(t)\subset F_{jq}$.
If we choose $\eta=\eta_{jq}$ sufficiently small, then $\Delta\cap F_{jq}(t)$ will be a Jordan arc, $\overline\Sigma_{3j+1}\cap F_{jq}\cap\partial\Delta=\emptyset$, and, by (6.1),
$$
x(s)\in\cases M(a',b'),&\text{if $s\in F_{jq}$}\\
M(c),&\text{if $s\in F_{jq}\cap\overline\Sigma_{3j+1}$}.\endcases\tag6.2
$$
Since by Proposition 6.2 $x(\Delta(\var))\subset x(\Delta(\delta_1))\subset M(a)$, (6.2) implies $F_{jq}\cap\Delta(\var)=\emptyset$.
We also record that $F_{jq}(0)=\{s_{3j+3}\}$.

By continuity, each $r\in Q$ has a neighborhood $U_r\subset Q$ such that for $q\in\overline U_r$
$$
K_\tau(q) (s)\in\cases M(a',b'),&\text{if $s\in F_{jr}$}\\
M(c),&\text{if $s\in F_{jr}\cap\overline\Sigma_{3j+1}$}.\endcases\tag6.3
$$
Let $\{U_r\colon r\in R\}$ be a finite subcover, where $R\subset Q$.
At the price of refining, we can assume that the $U_r$ come from a rectangular grid, by slightly enlarging the cells of the grid.
In this case we can partition $\{ U_r\colon r\in R\}$ in $2^d$ families $\{U_r\colon r\in R_j\}$, $j=0,\ldots,2^d-1$, each consisting of disjoint sets.
Choose continuous functions $\chi_r\colon Q\to [0,1]$, supported in $U_r$, so that $\sup_{r\in R}\chi_r\equiv 1$, and define for $(t,q)\in [\tau,\overline\tau]\times Q$ 
$$
G_{tq}=\Delta\backslash \bigcup\{F_{jr}\bigl((t-\tau)\chi_r (q)\bigr)
\colon q\in U_r,\ r\in R_j\}.\tag6.4
$$ 

Given $j$ (and $q$), there will be at most one $r\in R_j$ with $q\in U_r$, which ensures that the sets $F_{jr}((t-\tau)\chi_r (q))\subset\overline\Sigma_{3j+1}\cup\overline\Sigma_{3j+2}\backslash \{0\}$ in (6.4) are disjoint for $j=0,1,\ldots,2^d-1$.
It follows that $G_{tq}$ are Jordan domains, and their boundaries depend continuously on $t,q$.
(For fixed $q$, $G_{tq}$ can be described as $\Delta=G_{\tau q}$ minus a family of tentacles issued from finitely many points on $\partial\Delta$; with increasing $t$ the tentacles grow until, when $t=\overline\tau$, at least one of them becomes so large that it bridges some sector $\overline\Sigma_{3j+1}\cup\overline\Sigma_{3j+2}$.)
Note that $\Delta(\var)\subset G_{tq}$, since $F_{jr}\cap\Delta(\var)=\emptyset$.

Let $f_{tq}\colon\overline\Delta\to\overline G_{tq}$ be the homeomorphism that is holomorphic on $\Delta$ and satisfies $f_{tq} (0)=0$, $f'_{tq} (0)>0$.
By the Carath\'eodory--Rad\'o theorem [C, Ra], the map
$$
[\tau,\overline\tau]\times Q\times\overline\Delta \ni (t,q,s)\mapsto f_{tq} (s)\in\bC
$$
is continuous.
For $(t,q,s)$ as above, define
$$
K(t,q)(s)=\rho\{K_\tau (q) (f_{tq} (s))-h(\tau,q)+h(t,q)\}.\tag6.5
$$
The first term within the braces is contained in 
$K_\tau(q)(\overline\Delta)$, hence in some compact subset of $M(b')$ that 
is independent of $\overline\tau,t,q,s$; and the balance of the 
remaining terms is 
small if $\overline\tau$ is sufficiently close to $\tau$.
Hence $K(t,q)(s)\in M(b')$ is well defined, depends continuously on $(t,q,s)$ and holomorphically on $s\in\Delta$.
Similarly, (6.3) shows that for $q\in\overline U_r$
$$
K_\tau(q) (f_{tq}(s))\in\cases M(a',b'),&\text{if $s\in f^{-1}_{tq} (F_{jr}\cup\partial\Delta)$}\\
M(c),&\text{if $s\in f^{-1}_{tq}(F_{jr}\cap\overline\Sigma_{3j+1})$};\endcases
$$
in fact, in each case $K_\tau (q)(f_{tq}(s))$ is contained in a fixed compact subset of $M(a',b')$, resp.~$M(c)$, independent of 
$\overline\tau,t,q,r,s,j$.
Hence, again for $\overline\tau$ close to $\tau$,
$$
K(t,q)(s)\in\cases M(a',b'),&\text{if $s\in f^{-1}_{tq} (F_{jr}\cup\partial\Delta)$}\\
M(c),&\text{if $s\in f^{-1}_{tq} (F_{jr}\cap\overline\Sigma_{3j+1})$}.\endcases\tag6.6
$$
In particular, $K(t,q)(\partial\Delta)\subset M(a',b')$, since by (6.4)
$$
f_{tq}(\partial\Delta)=\partial G_{tq}\subset\partial\Delta\cup\bigcup^{2^d-1}_{j=0}\bigcup_{q\in U_r,r\in R_j} F_{jr}.
$$
As $K(t,q)(0)=h(t,q)$, $K$ is indeed a continuous lift of $h|[\tau,\overline\tau]\times Q$ to $D(a',b')$.

To finish the proof we consider $x=K(\overline\tau,q)$, and estimate the harmonic measure of $x^{-1}(M(c))\cap\partial\Delta$, using Lemma 5.3.
We have already noted that $G_{tq}\supset\Delta(\var)$.
Since $F_{jr}\subset\overline\Sigma_{3j+1}\cup\overline\Sigma_{3j+2}$ and 
$\overline\Sigma_{3j+1}\cap F_{jr}\cap\partial\Delta=\emptyset$, 
$G_{tq}$ also contains $\Sigma_{3j}$ and a one sided neighborhood of 
the arc of $\overline\Sigma_{3j+1}$, for every $j$.
Choose $j$ and $r\in R_j$ so that $\chi_r (q)=1$.
By (6.4) $\overline\Sigma_{3j+1}\cap\partial G_{\overline\tau q} = 
\overline\Sigma_{3j+1}\cap\partial F_{jr}$, and this latter contains the endpoint of $\Gamma_{jr}$, a point on the open segment $0s_{3j+1}$.
Thus we are in the situation of Lemma 5.3, and conclude that the harmonic measure of $\Sigma_{3j+1}\cap\partial G_{\overline\tau q}$, relative to $G_{\overline\tau q}$, is $>\delta_2$.
By invariance, the harmonic measure of 
$$
\Gamma=f_{\overline\tau q}^{-1} (\Sigma_{3j+1}\cap\partial G_{\overline\tau q})=f_{\overline\tau q}^{-1} (\Sigma_{3j+1}\cap \partial F_{jr}),
$$
relative to $\Delta$, is $>\delta_2$; but then we are done since (6.6) implies $x(\Gamma)\subset M(c)$.
\enddemo

\proclaim{Proposition 6.4}Assume that $h$ is smooth.
Given $a<a_0<b_0<b$ and $\var>0$, there is a $\delta>0$ with the following property.
Suppose with some $\tau\in [0,1-\delta]$ and $a<a'<a_0$, $b_0<b'<b$, $h(\tau,\cdot)$ has been lifted to $K_\tau\colon Q\to D(a',b')$.
If $\tau<\overline\tau\leq\tau+\delta$, then $h|[\tau,\overline\tau]\times Q$ has a lift $K\colon [\tau,\overline\tau]\times Q\to D(a'-\var,b'+\var)$ such that $K(\tau,\cdot)=K_\tau$ and $K(\overline\tau,\cdot)$ maps into $D'(a',b'+\var(\overline\tau-\tau))$.
\endproclaim

\demo{Proof}We can assume $\var < (b_0-a_0)/2$.
Pick $\delta_0 < (b_0-a_0)/4$ as in Lemma 4.3, $\delta_2$ as in Proposition 6.3, and choose constants $c_1,c_2 > 1$, $c_3,c_4,\delta>0$ so that 
$$
\align
& |h(t,q)-h(t',q)|\leq c_1 |t-t'|,\quad t,t'\in [0,1];\tag6.7\\
& |u(z)-u(\rho(w))|\leq c_2|z-w|,\quad z\in M(b),\ |z-w| < c_3;
\tag6.8\\
& c_1 c_2 e^{-c_4\delta_2} < \var\qquad\text{ and }\qquad\delta c_1 c_2 e^{c_4} < \min (\var,\delta_0,c_3).\tag6.9
\endalign
$$
Such a $\delta$ will do, as we now demonstrate.

By first lifting $h$ as in Proposition 6.3, but only over $[\tau,\tau']\times Q$ with some $\tau'\in(\tau,\overline\tau)$, we can reduce the proof to the situation where for every $r\in Q$ the initial lift $K_\tau(r)$ maps a set $\Gamma_r\subset\partial\Delta$ of harmonic measure $>\delta_2$ (relative to $\Delta$) into $M(a_0)$.
Choose complex polynomials $\theta_r\colon\bC\to\bC$ as in Proposition 5.2, i.e.
$$
\theta_r (0)=\delta_2,\qquad\text{ Re }\theta_r < \cases 1&\text{on $\overline\Delta$}\\ 0&\text{on $\partial\Delta\backslash\Gamma_r$}.\endcases
$$
Set $c=(a_0+b_0)/2$.
Each $r\in Q$ has a neighborhood $U_r\subset Q$ such that
$$
K_\tau (q) (\Gamma_r)\subset M(c)\quad\text{for}\quad q\in U_r.
$$
With a continuous partition of unity $\{\chi_r\}$ subordinate to $\{U_r\}$, $r\in Q$, let
$$
\theta(q,s)=\sum_{r\in Q}\chi_r (q)\theta_r (s),\quad (q,s)\in Q\times\oDelta.
$$
Thus $\theta(q,0)=\delta_2$, Re $\theta<1$ everywhere, and for 
$s\in\partial\Delta$, Re $\theta(q,s) < 0$ or $K_\tau(q)(s)\in M(c)$.

Now we apply Lemma 4.3 with $\kappa (q,s)=K_\tau (q)(s)$ and suitable 
$\alpha$, $\eta$ to obtain a continuous 
$\lambda\colon [\tau,\overline\tau]\times Q\times\oDelta\to M$, 
holomorphic on each slice $\{(t,q)\}\times\Delta$, such that 
$\lambda(\tau,\cdot)=K_\tau$, $\lambda(t,q,0)=K_\tau (q)(0)=h(\tau,q)$, 
and when $t\in[\tau,\overline\tau]$, $q\in Q$, $s\in\partial\Delta$ 
$$
\aligned
a'< u(\lambda (t,q,s)) < c+\delta_0,\quad & \text{if }K_\tau (q) (s)\in M(c);\\
a'+\delta_0<u(\lambda(t,q,s))<b',\quad & \text{if }K_\tau (q) (s)
\not\in M(c)\quad\text{or } t=\overline\tau.
\endaligned\tag6.10
$$
Define for $(t,q,s)\in [\tau,\overline\tau]\times Q\times\oDelta$
$$
K(t,q) (s)=\rho\bigl(\lambda (t,q,s)+
e^{c_4 (\theta(q,s)-\delta_2)} (h(t,q)-h(\tau,q))\bigr).\tag6.11
$$
Since by (6.7), (6.9),
$$
e^{c_4\theta(q,s)}|h(t)-h(\tau)|\leq e^{c_4\theta (q,s)} c_1 |t-\tau| < c_1 e^{c_4}\delta < c_3,\tag6.12
$$
(6.8) implies that $K$ is indeed well defined and holomorphic in $s\in\Delta$.
Also $K(t,q)(0)=h(t,q)$, so that $K$ lifts $h|[\tau,\overline\tau]\times Q$; and
$$
K(\tau,q)(s)=\lambda (\tau,q,s)=K_\tau (q)(s).
$$

To see where $K(t,q)$ maps $s\in\partial\Delta$, observe that by (6.11), (6.8), (6.12), (6.9)
$$
|u(\lambda(t,q,s))-u(K(t,q)(s))|\leq c_1 c_2 e^{c_4\theta(q,s)}|t-\tau|
<\var.\tag6.13
$$
Hence (6.10) implies that $K$ maps into $D(a'-\var,b'+\var)$, and that
$$ 
u(K(\overline\tau,q)(s))\geq u(\lambda(\overline\tau,q,s))-
\delta c_1 c_2 e^{c_4} > a'+\delta_0-\delta_0=a'.
$$
Further, in view of (6.13), (6.10), and (6.9) $u(K(\overline\tau,q)(s))$
can be estimated above by
$$
\cases u(\lambda(\overline\tau,q,s))+\delta c_1 c_2 e^{c_4}
< c+2\delta_0< b',&\text{if $K_\tau(q)(s)\in M(c)$}\\
u(\lambda(\overline\tau,q,s))+c_1 c_2 e^{-c_4\delta_2}(\overline\tau-\tau) < b'+\var(\overline\tau-\tau),&\text{if $K\tau(q)(s)\notin M(c)$}.
\endcases
$$
These estimates show that 
$K(\overline\tau,q)\in D(a',b'+\var(\overline\tau-\tau))$, and the 
proof is complete.
\enddemo

\demo{Proof of Theorem 4.1}First we show that given a smooth homotopy $h\colon [0,1]\times Q\to M(a)$ and a lift $H_0\colon Q\to D(a,b)$ of $h(0,\cdot)$, there is a lift $H\colon [0,1]\times Q\to D(a,b)$ of $h$ such that $H(0,\cdot)=H_0$.
Choose $\var > 0$ so that $H_0(Q)\subset D(a+3\var,b-3\var)$, set 
$$
a^*=a+2\var,\ a_0=a+3\var,\ b_0=b-3\var,\ b^*=b-2\var,
$$
and let $\delta=1/l$, $l\in\bN$, be as in Proposition 6.4.
We inductively construct lifts 
$$
H_j\colon [0,j/l]\times Q\to D(a,b),\quad j=0,1,\ldots,l
$$
of $h$, with the extra property that 
$$
H_j (j/l,q)\in D(a^*,b^*+j\var/l),\quad q\in Q.\tag6.14
$$
We are already given $H_0$.
Suppose $H_{j-1}$ has been constructed for $0< j\leq l$.
Proposition 6.4, applied with $\tau=(j-1)/l,\ \overline\tau=j/l,\ a'=a^*$, and $b'=b^*+(j-1) \var/l$, allows us to continue $H_{j-1}$ to a lift $H_j\colon [0,j/l]\times Q\to D(a,b)$, satisfying (6.14).
Proceeding in this manner we obtain the lift $H=H_l$ sought.

If $h$ is just continuous, we first choose $a< a^*<b^*<b$ so that 
$H_0(Q)\subset D(a^*,b^*)$.
We uniformly approximate $h$, $H_0$ by smooth maps 
$h'\colon [0,1]\times Q\to M(a^*)$, 
$H'_0\colon Q\to D(a^*,b^*)$ satisfying $\pi\circ H'_0=h'$.
By what we have proved, $h'$ has a continuous lift 
$H':[0,1]\times Q\to D(a^*,b^*)$ such that 
$H'(0,\cdot)=H'_0$.
We set
$$
\align
& H''(t,q)=H'(t,q)-H'_0 (q)+H_0(q)\text{ and}\\
& H(t,q)(s)=\rho\bigl(H''(t,q)(s)-H''(t,q)(0)+h(t,q)\bigr).
\endalign
$$
If the approximations were close enough, $H$ will be the lift we needed.
\enddemo

\subhead 7.\ Analytic Continuation in Mapping Spaces\endsubhead

\demo{Proof of Theorem 1.1}We take $M$ to be a closed submanifold of $\bC^n$ and $0\in M'$ the origin in $\bC^n$. Define a smooth strongly plurisubharmonic exhaustion function $u\colon M\to\bR$ by $u(z)=|z|^2$.
By the Docquier--Grauert theorem [DG] there are a pseudoconvex 
neighborhood $\Omega\subset\bC^n$ of $M$ and a holomorphic retraction 
$\rho\colon\Omega\to M$.
Thus, according to Proposition 2.1 $N=M_{S,A}\subset\Omega_{S,A}$ is a direct complex submanifold, in fact a holomorphic retract of 
$\Omega_{S,A}$; and $N'=M_{S,A}\subset N$ is an open subset.
We shall employ the method of sliding discs, Lemma 3.4, so that we have to 
consider a curve $[0,1]\ni t\mapsto x_t\in N$ starting at $x_0=0$, and show 
it has a regular lift $[0,1]\ni t\mapsto\xi_t\in N_\Delta$.
Choose real numbers $a<b$ so that
$$
\max_{z\in K} |z|^2 < a,
\qquad \max_{0\leq t\leq 1} \max_{s\in S} | x_t (s)|^2 < a,\tag7.1
$$
and $u$ has no critical value on $[a,b]$.
In particular $M[a,b]\subset M'$.

Suppose first $\theta\colon\Delta\to M$ is a proper holomorphic map 
such that $\theta(0)=0$ and $\theta(\Delta)\subset M'$.
For any $c\in (a,b)$ the level set $\{u\circ\theta < c\}\subset\Delta$ has simply connected components by the maximum principle.
We choose $c$ so that the boundaries of these components are smooth curves, and consider a biholomorphic map of $\Delta$ on the component of 0.
Composing $\theta$ with this map we obtain a proper holomorphic map 
$\Delta\to M(c)$, that extends to a smooth map 
$\theta_0\colon\oDelta\to M'$.
We can arrange that $\theta_0(0)=\theta(0)=0$.

We want to construct a continuous map $\Xi\colon [0,1]\times 
S\times\oDelta\to M$, holomorphic on slices $\{(t,s)\}\times\Delta$, 
that satisfies
$$
\align
\Xi(0,s,\sigma)\subset M',& \quad\text{if } 
(s,\sigma)\in S\times\oDelta,\tag7.2\\
\Xi(t,s,0)=x_t(s),&\quad\text{if } (t,s)\in [0,1]\times S,\tag7.3\\
\Xi([0,1]\times A\times\oDelta)=0,&\qquad \Xi([0,1]\times S\times\partial\Delta)\subset M'.\tag7.4
\endalign
$$
To this end we fix an $\var > 0$ so that if a $z\in M$ satisfies $|z|<\var$, then $\rho(z+\theta_0)$
defines a smooth map $\theta_z\colon\oDelta\to M(b)\cap M'$, holomorphic on $\Delta$.
Clearly $\theta_z(0)=z$.
There are a $\delta>0$ and a neighborhood $U\subset S$ of $A$ so that 
$|x_t (s)| < \var$ if $0\leq t\leq \delta$ or $s\in\overline U$.
We triangulate $S$ and choose $U$ to be the union of simplices 
of the triangulation.
With a smooth function 
$\chi\colon[0,1]\times S\to [0,1]$ satisfying
$$
\chi(t,s)=\cases 0&\text{if $t=0$ or $s\in A$}\\
1&\text{if $t\geq\delta$ and $s\notin U$},
\endcases
$$
we let, for $(t,s)\in ([0,1]\times\overline U)\cup ([0,\delta]\times S)$, 
$\sigma\in\oDelta$
$$
\Xi(t,s,\sigma)=\theta_{x_t(s)}(\chi(t,s)\sigma).
$$

The function $\Xi$ so far defined satisfies (7.2), (7.3), (7.4).
We now extend $\Xi$ to $[\delta,1]\times(S\setminus U)\times\oDelta$,
using Theorem 4.1. Let
$$
h(t,s)=x_t(s),\quad t\in [\delta,1],\quad s\in S\backslash U,
$$
a continuous map into $M(a)$.
We consider closed simplices $T\subset S\backslash U$ in the triangulation of $S$, and over $[\delta,1]\times T$ we lift the homotopy $h$ to a continuous $H_T\colon [\delta,1]\times T\to D(a,b)$, making sure that 
$H_T|T'=H_{T'}$ if $T'\subset T$.
We do this inductively on the dimension of $T$.
Suppose we have already lifted $h$ to all simplices of dimension 
$<\dim T=d$.
If $T\subset\partial U$, define
$$
H_T (t,s)=\Xi (t,s,\cdot),\qquad (t,s)\in [\delta,1]\times T.
$$
If $T\not\subset\partial U$, we let $Q=[0,1]^d$, and observe that the 
pair
$$
\bigl([\delta,1]\times T,\ ([\delta,1]\times\partial T)\cup 
(\{\delta\}\times T)\bigr)
$$
is homeomorphic to the pair $([0,1]\times Q,\ \{0\}\times Q)$.
The maps $H_{T'}$ for $T'\subset\partial T$ and $\Xi(\delta,\cdot,\cdot)$ 
together define a continuous map
$$
([\delta,1]\times\partial T)\cup (\{\delta\}\times T)\to D(a,b),
$$
which by Theorem 4.1 has a continuous extension
$$
H_T\colon [\delta,1]\times T\to D(a,b),
\qquad\text{with}\quad\pi\circ H_T=h|[\delta,1]\times T.
$$
Having constructed $H_T$ for all simplices in $S\backslash U$, we let
$$
\Xi (t,s,\cdot)=H_T (t,s),\qquad\text{if}\quad (t,s)\in[\delta,1]\times T;
$$
then (7.2), (7.3), (7.4) are satisfied.

If the mapping space $M_{S,A}$ was defined using
continuous maps $(S,A)\to (M,0)$, then the regular lift of $x_t$ sought is given by
$$
\xi_t(\sigma)=\Xi(t,\cdot,\sigma),\quad t\in [0,1],\quad \sigma\in\oDelta.
$$
To deal with spaces of $C^k$ or $W^{k,p}$ maps, we approximate 
$\Xi$ with a smooth map
$$
\Xi^*\colon [0,1]\times S\times\oDelta\to\bC^n,
$$
holomorphic on slices $\{(t,s)\}\times\Delta$, and vanishing on a
neighborhood of $[0,1]\times A\times\oDelta$.
If the approximation is close enough,
$$
\xi_t(\sigma)=\rho(\Xi^*(t,\cdot,\sigma)-\Xi^*(t,\cdot,0)+x_t)
$$
will define the required regular lift.
Theorem 1.1 now follows from Lemma 3.4.

On the other hand, if the map $\theta$ above does not exist, then 
$A=\emptyset$, and so
$M_{S,A}$ is independent of the choice of $0\in M'$.
We choose $a$ so that (in addition to earlier requirements, see (7.1)),
$K$ does not separate 0 from $\partial M(a)$.
Let $P\subset M(a)$ be the component of $0$.
Since the boundary of $P$ is strongly pseudoconvex, there is a proper holomorphic map $\theta\colon\Delta\to P$, mapping in a small neighborhood of a boundary point; in particular we can choose $\theta$ to avoid $K$.
Let a curve $[0,1]\ni t\mapsto z_t\in P\backslash K$ connect $z_0=\theta(0)$ with $z_1=0$, and define a curve $[0,2]\ni t\mapsto y_t\in P_{S,A}$ by
$$
y_t\cases\equiv z_t,&\text{if $0\leq t\leq 1$}\\
=x_{t-1},&\text{if $1\leq t\leq 2$}.\endcases
$$
We are in the situation of the first case of this proof, so $f$ has 
an analytic continuation along $t\mapsto y_t$; which then provides the 
analytic continuation along $t\mapsto x_t$ as well.
\enddemo

\demo{Proof of Theorem 1.2}Again we assume $M\subset\bC^n$ is a submanifold,
and $\Omega\subset\bC^n$ is a pseudoconvex neighborhood that
holomorphically retracts on $M$. By Proposition 2.1 
$M_{S,A}\subset\Omega_{S,A}$ is a direct complex submanifold, and
$\Omega_{S,A}$ holomorphically retracts on $M_{S,A}$. Furthermore,
$\Omega_{S,A}\subset\bC^n_{S,A}$ is pseudoconvex. Indeed, if $v$ is
a plurisubharmonic exhaustion function of $\Omega$, then
$$
w(x)=\sup_{s\in S}v(x(s)),\qquad x\in\Omega_{S,A},
$$
is plurisubharmonic and $w(x)\to\infty$ as $x$ tends to a point
in $\partial\Omega_{S,A}$.

First suppose that the mapping spaces are defined using Sobolev maps of class $W^{k,2}$.
Then $X=\bC^n_{S,A}$ is a Hilbert space, hence has an unconditional basis, and the theorem follows from Theorems 1.1 and 3.6.

Second, if $M_{S,A}, M'_{S,A}$ are built with another regularity class, permitted by Theorem 1.2, choose $k>(\dim S)/2$ so that $W^{k,2}$ contains this class.
We denote by $P$, $P'$ the versions of $M_{S,A}$, $M'_{S,A}$ built with 
$W^{k,2}$ maps. The inclusions
$$
M_{S,A}\hookrightarrow P,\qquad M'_{S,A}\hookrightarrow P'
$$
have dense image; also, they are homotopy equivalences by [P, Theorem 13.14], and so induce isomorphisms
$$
\pi_1(M_{S,A},x_0)\overset\approx\to\rightarrow \pi_1 (P,x_0),\qquad \pi_1(M'_{S,A},x_0)\overset\approx\to\rightarrow \pi_1 (P',x_0).
$$

We claim that the first of these isomorphisms maps 
$\Gamma=\Gamma (M_{S,A},M'_{S,A})$ into $\Gamma(P,P')$.
Indeed, let a closed curve $[0,1]\ni t\mapsto x_t\in M_{S,A}$ have 
homotopy class in $\Gamma$, and consider an $f\in\cO (P';\bC)$ together
with its analytic continuation $\{V_t,f_t\}$ along this curve.
The restriction $f|M'_{S,A}\in\cO(M'_{S,A};\bC)$ has analytic continuation 
$\{V_t\cap M_{S,A}, f_t|M_{S,A}\}$ along the same curve.
Since this latter continuation is single valued, $f_1(x)=f_0(x)$ for 
$x\in M_{S,A}$ near $x_0$; which, by density, implies $f_1(x)=f_0(x)$ for 
$x\in P$ near $x_0$.
Thus the analytic continuation of $f$ along $t\mapsto x_t$ is single 
valued, and so the homotopy class of this curve is in $\Gamma(P,P')$ as 
claimed.

In view of Proposition 3.5 (and Theorem 1.1), the homomorphisms 
considered therefore make up a commutative diagram
$$
\CD
\pi_1(M'_{S,A},x_0)@>\varphi>> \Gamma(M_{S,A},M'_{S,A})@>>> \pi_1(M_{S,A},x_0)\\
@V\approx VV @VVV @V\approx VV\\
\pi_1(P',x_0) @>\psi >> \Gamma(P,P') @>>> \pi_1 (P,x_0),
\endCD
$$
in which the horizontal arrows on the right are inclusions.
Since $\psi$ is onto by the first case of the proof, so is 
$\varphi$, q.e.d.
\enddemo

Now Gompf pointed out the following example.
If $S^2$ is embedded in $TS^2$ as the zero section, then
$$
0=\pi_2 (TS^2\backslash S^2)\to\pi_2 (TS^2)\approx\bZ\tag7.5
$$
is not onto.
This is relevant to Theorem 1.2 because $TS^2$ can be endowed with the structure of a Stein manifold $M$.
Letting $M'=TS^2\backslash S^2$, $S=S^1$, and $A\subset S^1$ a singleton, the homotopy groups in (7.5) are isomorphic to $\pi_1 (M'_{S,A})$, $\pi_1(M_{S,A})$, and so $\pi_1(M'_{S,A})\to\pi_1(M_{S,A})$ is not onto.
As a result, there is a holomorphic function on $M'_{S,A}$ that does not extend to a holomorphic function on $M_{S,A}$.
We conclude this paper by explicitly defining such a function.

Consider the quadric
$$
M=\{z\in\bC^3\colon\sum_{j=1}^3 z_j^2=1\}\quad\text{and}\quad K=M\cap\bR^3.
$$
It is known that 
$M'=M\setminus K$ is diffeomorphic to 
$W=(\bC^2\backslash \{0\})/w\sim\pm w$, hence 
$\pi_2(M')\approx\pi_2 (\bC^2\backslash \{0\})=0$.
The diffeomorphism $W\to M'$ is obtained by observing that the map
$$
\bC^2\backslash \{0\}\ni w
\mapsto (i(w_1^2+w_2^2),\ w_1^2-w_2^2,\ 2w_1 w_2)\in\bC^3
$$
induces a biholomorphism $W\to W_0=\{z\in\bC^3\backslash \{0\}\colon\sum_1^3 z_j^2=0\}$, and then composing this latter with the diffeomorphism
$$
W_0\ni z\mapsto (1+|\text{Re }z|^{-2})^{1/2}\text{ Re }z+i\text{ Im }z\in M'.
$$

We now take $S=\partial\Delta$ and for simplicity $A=\emptyset$.
We define the mapping spaces $M_{S,A}$, $M'_{S,A}$ using any regularity class $C^k\ (k\geq 1)$ or $W^{k,p}\ (k\geq 2)$; this ensures that elements of $M_{S,A}$ are $C^1$ maps.
The closed holomorphic 2--form
$$
\omega=z_1 dz_2\wedge dz_3-z_2 dz_1 \wedge dz_3+z_3 dz_1 \wedge dz_2|M
$$
induces a holomorphic function $f\colon M'_{S,A}\to\bC$ as follows.
Given $x\in M'_{S,A}$, i.e., a map $x\colon\partial\Delta\to M'$ homotopic to a constant, extend it to a $C^1$ map $\xi\colon\oDelta\to M'$, and put 
$$
f(x)=\int_\Delta\xi^*\omega.
$$
Since $\pi_2 (M')=0$, the integral above is independent of the choice of $\xi$ by Stokes' theorem. One checks that $f$ is holomorphic by
computing its differential. To this end, let $v\in T_xM'_{S,A}$. If, as
customary, we identify elements of $T_xM'_{S,A}$ with sections of the
induced bundle $x^*TM$, of the given regularity (see e.g., [L2, Proposition
2.2]), we find
$$
df(v)=\int_{\partial\Delta}\omega(v(s),dx(s)/ds)\,ds.\tag7.6
$$
This is complex linear in $v$, and $f$ is indeed holomorphic.
Now $f$ extends to a multivalued function $\tilde f$ on $M_{S,A}$, 
defined by the same recipe:\ if $x\in M_{S,A}$ has a 
$C^1$ extension $\xi\colon\oDelta\to M$, put 
$\tilde f(x)=\int_\Delta\xi^*\omega$.
However, this time the choice of $\xi$ matters, and in fact $\tilde f$ has no single valued branch.
Indeed, consider a smooth $\xi\colon\overline\Delta\to K=S^2$ 
that maps $\partial\Delta$ to a point, and is diffeomorphic on $\Delta$.
Let $\xi_t(\sigma)=\xi(1-t+\sigma t)$, $0\leq t\leq 1$, $\sigma\in\oDelta$, 
and $x_t=\xi_t|\partial\Delta$.
Thus $x_0=x_1$.
However, as $t$ increases from 0 to 1, the values $\tilde f(x_t)=\int_\Delta\xi_t^*\omega$ change from 0 to $\int_\Delta\xi_1^* \omega=\pm\int_K\omega$, and this latter is not zero since $\omega|K$ is an area form on $K$.
As the possible values of $\tilde f$ at any $x\in M_{S,A}$ differ by a period of $\omega$, we see that $\tilde f$ has no single valued branch along the curve $t\mapsto x_t$.

\enddocument